\documentclass{article}
\usepackage[margin=0.5in]{geometry}
\usepackage{amsmath}
\usepackage{amsthm}
\usepackage{amsfonts}
\usepackage{authblk}
\usepackage{titling}
\usepackage{graphicx}
\usepackage{indentfirst} 
\usepackage{hyperref}
\usepackage{bbm}
\usepackage{txfonts}
\usepackage[utf8]{inputenc}
\newtheorem{theorem}{Theorem}

\numberwithin{subsection}{section}
\numberwithin{theorem}{subsection}
\numberwithin{equation}{subsection}
\usepackage[T1]{fontenc}
\usepackage{titlesec}
\titleformat{\section}[hang]
{\bfseries\Large}{\thesection.}{0.5em}{}
\titleformat{\subsection}[runin] 
{\bfseries}{\thesubsection.}{0em}{}
\title{A Refinement of Carlson's Theorem}
\author{
Armen Vagharshakyan
}
\affil{
\small{Institute of Mathematics of Armenian Academy of Sciences}}
\date{}
\begin{document}
\maketitle
\begin{abstract}
    Carlson's theorem estimates the growth of an analytic function along the imaginary axis, provided that the function is zero at non-negative integers. We refine this theorem and describe not only the function's growth but also necessary and sufficient conditions in terms of its spectral measure. 
\end{abstract}
\section{Formulation}
\subsection{}
 Carlson's theorem (see \cite{C}) estimates the growth of an analytic function along the imaginary axis, provided that the function is zero at non-negative integers. Namely,
\begin{theorem}[Carlson]\label{tcarlson} Let $g$ be holomorphic in the closed right half plane,
\begin{equation}\label{chol}
    g\in Hol\left(Re(z)\geq 0\right),
\end{equation}
and let $g$ vanish at non-negative integers,
\begin{equation}\label{czeros}
    g(n)=0,\quad \text{for all } n\in \mathbb{Z}_+.
\end{equation}
Then the following dichotomy holds:
either $g\equiv 0$ or $g$ must grow substantially fast on the imaginary axis, 
\begin{equation}\label{cgrowth}
    \limsup_{y\in \mathbb{R},|y|\to+\infty}
    \frac{\ln|g(iy)|}{|y|}\geq \pi.
\end{equation}
\end{theorem}
In other words, under assumptions \eqref{chol} and \eqref{czeros}, Carlson's theorem provides a necessary condition \eqref{cgrowth}.
In this paper, under weaker assumptions \eqref{cholw} and \eqref{czerosw}, we find the necessary and sufficient condition \eqref{cnessuf}. That is, instead of estimating the function's growth, we describe its spectral measure. Namely, we prove:
\begin{theorem}\label{trefined}
Let the function $g$ be  holomorphic and of (at most) exponential type in the open right half plane, that is, for some $a\in \mathbb{R},b\geq 0$ we have
\begin{equation}\label{cholw}
    \sup_{x\geq \delta,y\in \mathbb{R}}\frac{\ln|g(x+iy)|}{ax+b|y|+\epsilon|z|}<+\infty,
    \quad \text{for all }\epsilon,\delta>0.
\end{equation}
Then the function $g$ vanishes at non-negative integers,
\begin{align}\label{czerosw}
    & g(n)=0,\quad \text{for all } n\in \mathbb{N},
    \\
    & \lim_{z\rightarrow 0}g(z)=0\nonumber
\end{align}
if and only if 
it can be represented as an integral of exponents via a spectral measure $\mu,$
\begin{equation*}
    g(z)=\int_{\mathbb{C}}e^{\omega z}d\mu(\omega)
\end{equation*}
with the additional restriction on $\mu$ that, if we introduce the circular-periodized equivalent of  $\mu$ by
\begin{align}\label{dnu}
& \nu\left(e^E\right)=\mu\left(E+\cup_{k \in \mathbb{Z}}\lbrace i\cdot2\pi  k\rbrace\right) ,\quad \text{for all } E\subset \mathbb{R}+i[-\pi,\pi)
\\
& \nu\left(\lbrace 0\rbrace\right)=0\nonumber,
\end{align}
then the measure $\nu$ has bounded support, and the Cauchy type integral $H$ of $\nu$'s sweeping measure belongs to the (complex-valued) Hardy space $\mathbb{H}_1$ and is zero at the origin, 
\begin{equation}\label{cnessuf}
   H\in \mathbb{H}_1,\quad H(0)=0.
\end{equation}
\end{theorem}
\subsection{}
The paper is composed as follows:
section \ref{shistory} contains a brief history of the subject,
section \ref{smeasures} is devoted to the proof of theorem \ref{trefined} by using Brown-Shields-Zeller theorem about balayage \ref{tbsz} and some basic facts about Hardy spaces in order to  specialize Morimoto's integral representation \ref{tmorimoto} to  assumption \eqref{czerosw},
section \ref{scarlson} contains a derivation Carlson's theorem from theorem \ref{trefined} by referring to Luzin-Privalov's boundary uniqueness principle \ref{tlp}, section \ref{spoly} discusses the analogue of our proof of Carlson's theorem for the polynomial case.
\section{Brief history}\label{shistory}
Carlson's theorem was first proved in \cite{C}. 
See \cite{H}, p.330 for an alternate proof by a certain summation formula. See \cite{YS2} for a proof by Plana's summation formula that uses Morimoto's integral representation indirectly.
\\\indent
An active area of research constitutes substituting the condition \eqref{czeros} to the assumption that $g$ vanishes on a different set. See e.g. \cite{GR}, \cite{F}, \cite{G}.
A general condition for when $g$ vanishes on a sequence is provided in \cite{MR}.
\\\indent
 Some of the applications of Carlson's theorem are:
the uniqueness of analytic continuation, 
Dyson conjecture in statistical mechanics, 
calculation of Selberg integral (as exposed in \cite{Y}).
It is also closely related to completeness of systems of exponentials (see \cite{Kh}).
\section{Proof of the Main Theorem  \ref{trefined}}\label{smeasures}
\subsection{}\label{ssmorimoto}
We define the classes $Exp_{a,b}$ of exponential functions in the right half plane and then describe a Paley-Wiener type theorem for them. 
For $-\infty<a<+\infty, \; 0\leq b<+\infty,$ define $Exp_{a,b}$ to be the space of functions $g$ holomorphic in the open right half plane, whose growth is bounded by the constants $a,b$ as follows:
\begin{equation*}
    Exp_{a,b}=\left\lbrace
    g(x+iy)\in Hol(x>0)\colon
      \forall\epsilon,\delta>0
    \left[
    \sup_{x\geq \delta,y\in \mathbb{R}}\frac{\ln|g(x+iy)|}{ax+b|y|+\epsilon|z|}<+\infty
   \right]
   \right\rbrace
\end{equation*}
\begin{theorem}[Morimoto]\label{tmorimoto}
For any $g\in Exp_{a,b}$ and $\epsilon>0$ there exists a finite (complex-valued) spectral measure $\mu\colon \mathbb{C}\rightarrow \mathbb{C}$ with (possibly unbounded) support
\begin{equation}\label{cmubound}
supp\left(\mu\right)\subset (-\infty,a+\epsilon]+i[-b-\epsilon,b+\epsilon]
\end{equation}
such that the function $g$ may be expressed as the following integral of exponents
\begin{equation}\label{dmu}
    g(z)=\int_{\omega\in\mathbb{C}} e^{\omega z}d\mu(\omega).
\end{equation}
\end{theorem}
\begin{proof}
For a proof, see theorem 5.1 in \cite{M} together with proposition 2.1.1 in \cite{SM}. Also, to recover the exact phrasing that we used, insert 
\begin{equation*}
k^{\prime}=0,\quad\epsilon^{\prime}=1,\quad e^{\frac{\epsilon^{\prime}}{2}|Re(\omega)|}d\mu(\omega)=d\mu(\omega)
\end{equation*}
in theorem 1.1 and proposition 2.4 in \cite{YS1}.
\end{proof}
\subsection{}
While keeping the notations introduced in \ref{ssmorimoto},  define the measure $\nu$ as the circular-periodized equivalent of the measure $\mu.$ Namely,
relate the measures $\mu$ and $\nu$ by \eqref{dnu}.
We now translate the properties of $\mu$ into those of $\nu.$ Firstly, it follows from condition \eqref{cmubound} that the measure $\nu$ has bounded support.
Also, in view of integral representation \eqref{dmu} of the function $g$  and definition \eqref{dnu} of the measure  $\nu,$ the conditions \eqref{czerosw} are equivalent to the following condition on the measure $\nu:$
\begin{equation}\label{cnu}
   \int_{\mathbb{C}}\zeta^n d\nu=0,\quad \text{for all }n\in\mathbb{Z}_+.
\end{equation}
\subsection{}
Assume that the measure $\nu$ is non-trivial,
$
    \nu\not\equiv 0.
$
Then the number 
\begin{equation}\label{er}
    r=\sup_{z\in supp(\nu)}|z|
\end{equation}
is well-defined. The measure $\nu$ has bounded support, so that
$r$ is finite.
Consider the open disc
\begin{equation}\label{ed}
D=\lbrace \zeta\in \mathbb{C}\colon |\zeta|< r\rbrace.
\end{equation}  
We will use a ``balayage'' type theorem \ref{tbsz} to sweep the measure $\nu$ (whose support lies in the closed disc $\bar{D}$) into the boundary $\partial{D}$ of that disc. By sweeping, we will reduce the case of the planar measure $\nu$ to the case of a linear measure. Indeed, introduce two measures associated with the measure $\nu:$ the internal measure $\nu_{int}$ - the planar measure which is the restriction of $\nu$ to $D,$ and the external measure $\nu_{ext}$ - the linear measure which is induced by $\nu$ on $\partial{D}.$ 
\begin{theorem}[Brown, Shields, Zeller]\label{tbsz}
For a bounded measure $\nu_{int}$ on the open disc $D$,
there exists a function $h_{int}\in L_1\left(\partial{D}\right)$ such that for all bounded analytic functions $f$ on $D$ we have
\begin{equation}\label{ebalayage}
    \int_{D}fd\nu_{int}=\int_{\partial{D}}f \cdot h_{int},
\end{equation}
\end{theorem}
\begin{proof}
For a proof see \cite{BSZ}. Also, to recover the exact phrasing that we used, insert 
$\mu=\nu$
into formula (4.23.1) of \cite{RS}.
\end{proof}
\subsection{}\label{sspoisson}
As a consequence of theorem \ref{tbsz}, we may write
\begin{equation}\label{epreriesz}
0
=
\int_{\mathbb{C}}\zeta^n d\nu
=
\int_{\bar{D}}\zeta^n d\nu
=
\int_{\partial{D}}\zeta^n \cdot\left(h_{int}+d\nu_{ext}\right)
,
\quad 
\text{for all } n\in \mathbb{Z}_+,
\end{equation}
where the first equality is condition \eqref{cnu},
the second equality is implied from the definition of $D$ (see \eqref{ed},\eqref{er}),
and the third equality comes from the sweeping argument \eqref{ebalayage}.
By the F. and M. Riesz theorem (see e.g. \cite{O}),  equality \eqref{epreriesz} implies that the measure $\nu_{ext}$ is absolutely continuous on $\partial {D}.$ Denote
its Radon-Nikodym derivative (with respect to the linear Lebesgue measure on $\partial{D}$) by $h_{ext}\in L_1(\partial{D}).$ For brevity, denote
$h=h_{int}+h_{ext}\in L_1(\partial{D})$ 
and rewrite
\eqref{epreriesz} as
\begin{equation}\label{eh}
\int_{\partial{D}}\zeta^n \cdot h
=0
,
\quad 
\text{for all } n\in \mathbb{Z}_+.
\end{equation}
Denote by $H$ the Cauchy type integral of the function  $h,$ or equivalently
\begin{equation*}
    H(\zeta)=\sum_{n=0}^{+\infty}a_n \zeta^n, \quad a_n=\fint_{\partial{D}} \zeta^{-n}\cdot h.
\end{equation*}
By theorem 3.12 in \cite{Ka}, the condition \eqref{eh} is equivalent to claiming that $H$ is in the (complex-valued) Hardy space $\mathbb{H}_1$ and is zero at the origin, that is condition \eqref{cnessuf} holds.
Moreover, a review of section \ref{smeasures} shows that conditions \eqref{czerosw} imposed on the function $g$ are in fact equivalent to the condition \eqref{cnessuf}.
\section{Derivation of Carlson's theorem}\label{scarlson}
\subsection{}\label{ssboas}
Let $g$ be holomorphic in the closed right half plane, $g\in Hol(Re(z)\geq 0)$.
Recall that the indicator of $g$ is defined to be the following function:
\begin{equation*}
    h_g(\theta)=\limsup_{r\rightarrow +\infty}\frac{\ln\left|g\left(re^{i\theta}\right)\right|}{r},\quad -\pi/2\leq \theta\leq \pi/2.
\end{equation*}
For brevity denote
$I^{\ast}=h_g(\pi/2),I_{\ast}=h_g(-\pi/2).$
Assume that 
$I^{\ast},I_{\ast}<+\infty.$
From trigonometric convexity of the indicator, it follows that for some $a\in \mathbb{R}$ we have
$g\in Exp_{a,\max\left(I^{\ast},I_{\ast},0\right)}.$
Pick any $\epsilon>0.$ Let $\mu$ be the corresponding spectral measure provided by theorem \ref{tmorimoto}. Denote the effective endpoints of the integral in formula \eqref{dmu} to be
\begin{equation*}
    b^{\ast}=\sup_{\omega\in supp(\mu)}Im(\omega),
\quad
    b_{\ast}=\inf_{\omega\in supp(\mu)}Im(\omega).
\end{equation*}
From the description \eqref{cmubound} of $\mu$'s support, we estimate
\begin{equation}\label{ebleqi}
\max\left(\left|b_{\ast}\right|,\left|b^{\ast}\right|\right)
\leq     
\max\left(I^{\ast},I_{\ast},0\right)+\epsilon.
\end{equation}
\subsection{}\label{snewproof}
We now deduce Carlson's theorem \ref{tcarlson} from theorem
\ref{trefined}. Indeed, let the inequality \eqref{cgrowth} claimed by Carlson's theorem not hold, that is
\begin{equation}\label{cindicator}
    \max\left(I_{\ast},I^{\ast}\right)<\pi.
\end{equation}
By  \eqref{ebleqi}, the restriction \eqref{cindicator} implies
\begin{equation}\label{cond_spectrum}
    -\pi<b_{\ast}\leq b^{\ast}<\pi.
\end{equation}
We now reduce the case of two-dimensional spectral measure $\mu$ to that of one-dimensional, as claimed in \eqref{egtohh}. Specifically, assume $g\not\equiv 0,$ so that the disc $D$ is well-defined by \eqref{ed}. We may write
\begin{equation}\label{egtoh}
    g(z)=
    \int_{\mathbb{C}}e^{\omega z}d\mu(\omega)=
    \int_{\bar{D}}\zeta^z d\nu(\zeta)=
    \int_{\partial{D}}\zeta^z \cdot h(\zeta)=
    \int_{\ln(r)+i[-\pi,\pi)}e^{\omega z} \cdot h\left(e^{\omega}\right),\quad Re(z)>0.
\end{equation}
Here, the first equality is the integral representation \eqref{dmu}.
The second equality doesn't hold for arbitrary measures $\mu,$ as the change of variables involved in passing from integration over $\mu$ to integration over $\nu$ may not be one-to-one, or, in other words, the function $\zeta\rightarrow \zeta^z$ is, generally speaking, multivariate. However, due to restriction \eqref{cond_spectrum} on the support of measure $\mu$, the mapping $\zeta=e^{\omega}$ is one-to-one on $\omega\in supp(\mu).$ Also, due to restriction \eqref{cond_spectrum}, while integrating a function over the measure $\nu,$ the values of that function on the ray
$\arg(\zeta)=-\pi$ are not taken into account. Thus, while talking about the integral of the function $\zeta\rightarrow\zeta^z$ over the measure $\nu,$ we can think of the branch of that function, generated by splitting the complex plane by the ray $\arg(\zeta)=-\pi.$  
The third equality follows from applying the theorem \ref{tbsz} about balayage to the bounded analytic function $\zeta\rightarrow \zeta^z$; here,
the fact that it is indeed a bounded analytic function follows from the restriction $Re(z)>0.$
The last equality is obtained by a change of variable.
For brevity, rewrite \eqref{egtoh} as an integral representation of the function $g$ by a one-dimensional spectral measure, 
\begin{equation}\label{egtohh}
     g(z)=\int_{\ln(r)+i[-\pi,\pi)}e^{\omega z} \cdot    h\left(e^{\omega}\right),\quad Re(z)>0.
\end{equation}
By \cite{T}, p.323 (or by corollary 6.9.4 in \cite{B}, p.108), the restriction \eqref{cindicator} implies
\begin{equation}\label{cprelp}
    -\pi
    <\inf_{\zeta\in   supp(h)}arg(\zeta)
    \leq 
    \sup_{\zeta\in supp(h)}arg(\zeta) 
    <\pi
\end{equation}
Recall Luzin-Privalov's boundary uniqueness principle for the (complex-valued) Hardy space $\mathbb{H}_1$ (see e.g. \cite{DT}, pp. 73-74, or \cite{Ka}, p.102 for our exact formulation):
\begin{theorem}[Luzin-Privalov]\label{tlp}
Let $f\in \mathbb{H}_1.$ Then the following dichotomy holds: either $f\equiv 0$ or $f$ can vanish only on a subset of the unit circle $|z|=1$ of measure zero.
\end{theorem}
Remarkably, Luzin-Privalov's principle is phrased as a dichotomy, just as Carlson's theorem does.
By   theorem \ref{tlp}, the restriction \eqref{cprelp} implies $h\equiv 0$. By \eqref{egtohh}, we have $g\equiv 0,$ yet we assumed the converse.
The proof of Carlson's theorem is now complete.
\section{Relation to the polynomial case}\label{spoly}
\subsection{}\label{sspoly}
We remark that if $p$ is a certain polynomial with (at least) $n$ zeros, then the following dichotomy holds: either $p\equiv 0$ or $p$ must grow substantially fast, namely
\begin{equation*}
    \lim_{x\to\infty}\frac{\ln|p(x)|}{\ln|x|}\geq n.
\end{equation*}
Carlson's theorem is the analogue of this remark
for analytic functions.
\subsection{}
Remark \ref{sspoly} about polynomials may be proved by polynomial remainder (little Bezout) theorem. Similarly, Carlson's theorem may be proved by using the condition \eqref{czeros} to express $g$ as an infinite product, and by estimating the growth of $g$ from this expression (see e.g. \cite{H}). 
\\\indent
However, there is another way to prove remark \ref{sspoly} about polynomials.
Namely, the polynomial $p$ is a sum of monomials; further, the rate of growth of that sum coincides with that of its main term; finally, the proof will follow 
from the maximum principle - if the polynomial $p$ is non-zero, then the polynomial $p$ and its main term have the same number zeros in $\mathbb{C}$. Hence, remark \ref{sspoly} follows.
In \ref{snewproof} we proved Carlson's theorem  in a similar way. 
Namely, the holomorphic function $g$ turned out to be an integral of exponents; further, the rate of growth of that integral was determined by its effective endpoints; finally,
the proof followed from Luzin-Privalov's boundary uniqueness principle, an analogue of maximum principle.
\subsection{}
We also note that, if (some of the) zeros of a non-zero polynomial are given explicitly, say,
\begin{equation}\label{epzeros}
    p(0)=p(1)=\dots=p(n-1)=0,
\end{equation}
then not only we can estimate $p$'s growth from below, but additionally, by invoking Vieta's formulas, we can 
claim that certain relations between coefficients of $p$'s expansion into monomials are equivalent to \eqref{epzeros}. Our main theorem \ref{trefined}, being a refinement of Carlson's theorem, may be viewed as a generalization of this fact for analytic functions. For example, \ref{ssdiscrete} claims certain relations between coefficients of an analytic function's expansion into exponents.
\subsection{}\label{ssdiscrete}Note that the only discrete measure, that is also absolutely continuous, is the zero measure.
Hence, by specializing theorem \ref{trefined} to discrete measures, we obtain the following remark: a finite sum of exponents
\begin{equation*}
    g(z)=\sum_{k=1}^N c_k e^{i\omega_k z},\quad  \lbrace\omega_k\rbrace_{k=1}^N\subset R
\end{equation*}
equals zero at non-negative integers if and only if
\begin{equation*}
  \sum_{k\colon \omega_k\; mod\; 2\pi=\omega } c_k=0, \quad \text{for all } \omega\in [0,2\pi).
\end{equation*}


\begin{thebibliography}{Xyz12}
\bibitem{C}
F. Carlson, 
``Sur une classe de series de Taylor,” 
Thesis, Uppsala, Sweden (1914).
\bibitem{H}
G. Hardy,
``On two theorems of F. Carlson and S. Wigert," 
Acta Mathematica, vol. 42 (1920), pp. 327–339. 
\bibitem{YS2}
K. Yoshino, M. Suwa,
``A new proof of Carlson's theorem by Plana's summation formula,"
Surikaisekikenkyusho Kokyuroku (2001), pp. 166-174.
\bibitem{GR}
R. Gervais, Q. Rahman, 
``An extension of Carlson's theorem for entire functions of exponential type", 
Trans. Amer. Math. Soc., vol. 235 (1978), pp. 387–394.
\bibitem{F}
W. Fuchs,
``A Generalization of Carlson's Theorem,"
J. of the Lond. Math. Soci., vol. 21, no. 2 (1946) pp. 106–110.
\bibitem{G}
A. Goldberg,
``The integral over a semi-additive measure and its application to the theory of entire functions. I,"
Math. Sb., vol 58, no. 3 (1962), pp. 289-334
\bibitem{MR}
P. Malliavin, L.A. Rubel,
``On small entire functions of exponential type with given zeros," 
Bulletin of the Mathematical Society of France, Tome 89 (1961), pp. 175-206. 
\bibitem{Y}
K. Yoshino
``On Carlson’s theorem for holomorphic functions,"
Algebraic Analysis (Edited by M.Kashiwara and T.Kawai), Vo1.2, Academic Press (1988), pp. 943-950.
\bibitem{Kh}
B. Khabibullin,
``Completeness of systems of exponents and uniqueness sets," (in Russian) 3rd ed.,
Ufa (2011), ISBN 978-5-7477-2540-9.
\bibitem{DT}
E. Dolzhenko, G. Tumarkin,
``N. N. Luzin and the theory of boundary properties of analytic functions,"
Uspekhi Mat. Nauk, vol. 40, issue 3/243/ (1985), pp. 71-84.
\bibitem{Ka}
Y. Katznelson,
``An Introduction to Harmonic Analysis," 3rd corr. ed., 
Cambridge University Press (2002).
\bibitem{M}
M. Morimoto, 
``Analytic functionals with non-compact carrier", 
Tokyo J.Math., vol.1 (1978), pp. 77-103.
\bibitem{SM}
P. Sargos, M. Morimoto, 
``Transformation des fonctionelles analytiques a porteurs non compacts," 
Tokyo J. Math., vol. 4 (1981), pp. 457-492.
\bibitem{YS1}
K. Yoshino, M. Suwa,
``Plana's summation formula for holomorphic functions of exponential type,"
Suuri Kaiseki Kenkyuujo koukyuuroku (2000), pp. 180-189.
\bibitem{BSZ}
L. Brown, A. Shields and K. Zeller,
``On absolutely convergent exponential sums,"
Trans. Amer. Math. Soc., vol. 96 (1960), pp. 162-183.
\bibitem{RS}
L. Rubel, A. Shields, 
``The space of bounded analytic functions on a region,"
Ann. Ins. Fourier (Grenoble), vol. 16 (1966), pp. 235-277.
\bibitem{O}
B. Oksendal, 
``A short proof of the F. and M. Riesz theorem,"
Proc. AMS, vol. 30, no. 1 (1971), p. 204.
\bibitem{T}
E. Titchmarsh,
``Introducton to the Theory of Fourier Integrals"
Oxford university Press, 1937.
\bibitem{B}
R. Boas, 
``Entire Functions,” 
Academic Press, New York, 1954.
\end{thebibliography}
\end{document}